\documentclass{article}
\usepackage[a4paper, margin=3cm]{geometry}
\usepackage[utf8]{inputenc}
\usepackage[english]{babel}

\usepackage{soul}
\usepackage{amsthm}
\usepackage{amssymb}
\usepackage{indentfirst}
\usepackage{amsmath}
\usepackage{graphicx, color}
\usepackage{enumerate}
\usepackage{algorithm}
\usepackage{algpseudocode}
\usepackage{authblk}
\usepackage{xcolor}
\usepackage{colortbl}
\usepackage{todonotes}{}

\usepackage{comment}

\usepackage{tikz}

\newtheorem{theorem}{Theorem}[section]
\newtheorem{lemma}[theorem]{Lemma}

\newtheorem{observation}[theorem]{Observation}
\newtheorem{conjecture}[theorem]{Conjecture}
\newtheorem{claim}[theorem]{Claim}

\newcommand{\zig}{\operatorname{zig}}
\usepackage{ulem}

\newcommand{\piso}[1]{\left\lfloor#1\right\rfloor}
\newcommand{\conjunto}[1]{\left\{#1\right\}}
\newcommand{\z}{\mathcal{Z}}
\newcommand{\C}{\mathrm{\mathcal{C}}}

\floatname{algorithm}{Procedure}

\newcommand{\qitem}[1]{\noindent\leavevmode\hangindent1.5\parindent%
  \noindent\hbox to1.5\parindent{#1\hss}\ignorespaces}

\newenvironment{claimproof}{ \trivlist
	\item[\hskip\labelsep
	\textit{Proof of the claim}.]\ignorespaces
}{\hfill$\vartriangleleft$\medskip
	
}

\title{Totally odd subdivisions in Kneser graphs}

\author[1]{Henry Echeverr\'ia}
\author[1,2]{Andrea Jim\'enez}
\author[3]{Suchismita Mishra}
\author[4,5]{Adri\'an Pastine}
\author[1]{Daniel~A.~Quiroz}
\author[1]{Mauricio Y\'epez}

\affil[1]{\small Instituto de Ingenier\'ia Matem\'atica-CIMFAV, Universidad de Valpara\'iso, Chile.} 

\affil[2]{Millennium Nucleus for Social Data Science (SODAS), Santiago, Chile}

\affil[3]{Institute of Mathematical Science, Chennai, India}

\affil[4]{Instituto de Matem\'atica Aplicada San Luis (UNSL-CONICET), Argentina.}

\affil[5]{Departamento de Matem\'atica, Universidad Nacional de San Luis, Argentina.}

\date{}

\begin{document}
	
\maketitle
	
\begin{abstract}

As evidence for the Odd Hadwiger Conjecture, Simonyi and Zsbán (2010) showed that every Kneser graph $G$ with large enough order (compared to $\chi(G)$) contains a totally odd subdivision of $K_{\chi(G)}$. A recent result of Steiner (2024), shows that every Schriver graph, and thus every Kneser graph, satisfies the Odd Hadwiger Conjecture, that is, it contains $K_{\chi(G)}$ as an odd minor. We strengthen these results for Kneser graphs in two ways. We show that for every $t\ge 8$, there are $t$-chromatic Kneser graphs that contain arbitrarily large complete totally odd subdivisions (and thus, odd minors). We also show that every Kneser graph contains a totally odd subdivision of $K_{\chi(G)}$. 

Kneser graphs are the prime example of graphs having chromatic number equal to its topological lower bounds. Motivated by our main results, we also study totally odd immersions on graphs with this property, proving, in particular, that if the chromatic number of $G$ is equal to any of its topological lower bounds, then $G$ contains a totally odd immersion of $K_{\lfloor \chi(G)/2 \rfloor +1}$. This gives evidence for the immersion-analogue of the Odd Hadwiger Conjecture.

\end{abstract}


\section{Introduction}

\noindent
All graphs in this paper are finite, undirected and loopless.

Hadwiger's conjecture \cite{Hadwiger} is an ambitious conjecture that aims to greatly generalise the Four Colour Theorem; it states that every graph $G$ contains $K_{\chi(G)}$ as a minor. This famous conjecture is trivial for $\chi(G)\le 3$, was verified for $\chi(G)=4$ by Hadwiger~\cite{Hadwiger} and Dirac~\cite{D52}, and cases $\chi(G)=5,6$ where reduced to the Four Colour Theorem~\cite{Wagner,H6}. For large values of $\chi(G)$, Kostochka~\cite{Kostochka} and Thomason~\cite{Thomason} showed that any graph $G$ with no $K_t$ minor is $O(t\sqrt{\log t})$-degenerate (which is best possible) and thus that $\chi(G)\in O(t\sqrt{\log t})$. These asymptotics were improved by Norin, Postle and Song~\cite{breaking} and later by Delcourt and Postle~\cite{DP21} who showed that every graph $G$ with no $K_t$ minor  satisfies $\chi(G)\in O(t\log{\log t})$.

A strengthening of Hadwiger's conjecture was proposed by Gerards and Seymour in terms of odd minors (see~\cite{JT95}). The \underline{odd minor} relation is a restriction of the minor relation that preserves parities of cycles: it allows for the operations of deletion of vertices and edges, and of contraction of edge-cuts. The strengthened conjecture is known as the Odd Hadwiger Conjecture, and reads as follows.

\begin{conjecture}[Gerards and Seymour]\label{conj:oddhad}
    Every graph $G$ contains $K_{\chi(G)}$ as an odd minor.
\end{conjecture}

This conjecture is also trivial for the case $\chi(G)\le 3$, and was verified for $\chi(G)=4$ by Catlin~\cite{C79}. In terms of asymptotics, after various bounds were given (e.g.~\cite{GGRSV06,NorinSong}), Steiner~\cite{S22} recently tied the asymptotics of the Odd Hadwiger Conjecture to those of Hadwiger's conjecture. 

Being such a hard conjecture, it is interesting to verify Conjecture~\ref{conj:oddhad} for specific graph classes. In this sense Simonyi and Zsb\'an~\cite{SimonyiZsban} showed that every Schrijver graph\footnote{Set $[n]=\{1,2,\dots ,n\}$. For integers $n,k$ with $n\geq k$, the \underline{Kneser graph} $KG(n,k)$ is the graph having as vertex set all the subsets of $[n]$ of size $k$, and having an edge between two vertices if their intersection is empty. The \underline{Schrijver graph} $SG(n,k)$ is the subgraph of $KG(n,k)$ induced by all $k$-subsets of $[n]$ having no pair of consecutive elements, where 1 and $n$ are also considered consecutive. As proved by Schrijver~\cite{Schrijver} the chromatic number of $SG(n,k)$ is equal to that of $KG(n,k)$ and 
deleting any of its vertices reduces the chromatic number.} $G$ (and thus every Kneser graph) with large enough order compared to $\chi(G)$ satisfies Conjecture~\ref{conj:oddhad}. In fact, they showed something stronger: every such graph $G$ contains a totally odd subdivision of $K_{\chi(G)}$ (a \underline{totally odd subdivision} is a subdivision where every path joining two terminals is of odd length). Working directly with odd minors, Steiner~\cite{Steiner}, was able to show that every Schrijver graph satisfies Conjecture~\ref{conj:oddhad}. In this paper, we give two results which strengthen all of this in the case of Kneser graphs; these are the following.

\begin{theorem}\label{theo:Kneser}
For every $s\ge 8$ and $t\ge 1$, there exist a Kneser graph $G$ with $\chi(G)=s$ which contains $K_t$ as a totally odd subdivision.
\end{theorem}

\begin{theorem}\label{theo:Kneser2}
Every Kneser graph $G$ contains $K_{\chi(G)}$ as a totally odd subdivision.
\end{theorem}

Our next result studies substructures similar to totally odd subdivisions in graphs having their chromatic number equal to any of it topological lower bounds, of which Kneser graphs are prime examples.

Lov\'asz proof of Kneser's conjecture \cite{Lovasz}, which gives a tight lower bound on the chromatic number of Kneser graphs, was groundbreaking in its use of topological methods for such a combinatorial problem. 
Since Lovász's work, many topological parameters have been introduced as lower bounds for the chromatic number. In \cite{MZ}, Matou\v{s}ek and Ziegler made a hierarchy of these parameters, which was recently updated by Daneshpajouh and Meunier \cite{DaneshpajouhAndMeunier}. Out of these graph parameters, one introduced by Simonyi and Tardos~\cite{SimonyiTardos},  the 
zigzag number, denoted $\zig(G)$, stands out for two reasons. The first reason is its position in the hierarchy mentioned above: $\zig(G)$ is an upper bound to all other topological bounds of~$G$, and still a lower bound for $\chi(G)$, therefore proving that something holds for graphs $G$ with $\zig(G)=\chi(G)$, implies proving that it holds for graphs where \textit{any} topological bound meets the chromatic number. The second reason why the zigzag number stands out is that, while it arose from the study of topological bounds, it is a purely combinatorial parameter. Moreover, the fact that it is a lower bound for the chromatic number, does not depend on topological arguments, but can be easily deduced from its definition, which we now give. 

 Let $c$ be a proper colouring of a graph $G$. A \emph{zigzag} in $(G,c)$ is a  sequence of vertices $z_1,z_2,z_3, \dots z_t$ satisfying the following.
    
    \begin{enumerate}[(i)]
       
        \item $c(z_i) < c(z_{i+1})$, for all $1 \leq i <t$.
        \item For each even number $i\in [t]$, 
 $z_i$ is adjacent to  all vertices $z_j$  with $j \in [t]$ an odd number.    \end{enumerate} 
    
    The \emph{zigzag number of $(G,c)$}, denoted $\zig(G,c)$, is defined as the size of a maximum zigzag of~$(G,c)$. While the \emph{zigzag number of $G$},  $\zig(G)$, is the minimum $\zig(G,c)$ over all the proper colourings $c$ of $G$.

 After proving that large enough Shrijver graphs satisfy the Odd Hadwiger Conjecture, Simonyi and Zsb\'an~\cite{SimonyiZsban} asked whether a more general result could be obtained in terms of the zigzag number. They asked if there is a positive constant $\varepsilon$ such that every graph~$G$ with $\zig(G)\ge t$ contains an odd minor of $K_{\lceil \varepsilon t\rceil}$. This was recently answered in the positive by Steiner~\cite{Steiner} who showed that it holds with constant $\varepsilon=\frac 12$. It is unlikely that such a result could be strengthened by replacing odd minors with totally odd subdivisions, but we wonder if an analogous result could be obtained for totally odd immersions\footnote{A graph~$G$ contains another graph $H$ as an \emph{immersion} if there exists an injection $\varphi\colon V(H)\rightarrow V(G)$ such that:
\begin{itemize}
\item For every $uv\in E(H)$, there is a path in $G$, denoted $P_{uv}$, with endpoints $\varphi(u)$ and~$\varphi(v)$.
\item The paths in $\{P_{uv} \mid uv\in E(H) \}$ are pairwise edge disjoint.
\item The vertices of $\varphi(V(H))$, called \emph{terminals}, do not appear as internal vertices on paths~$P_{uv}$.
 \end{itemize}
If the paths are internally vertex-disjoint, then $G$ contains $H$ as a \emph{subdivision}.  Thus, if $G$ contains $H$ as a subdivision, then it contains $H$ as an immersion. If $G$ contains $H$ as an immersion (resp. subdivision) where all the corresponding paths are of odd length, then we say that $G$ contains $H$ as a \emph{totally odd immersion} (resp. \emph{totally odd subdivision}).}. Since minors and immersion are incomparable, such a result would be incomparable to that of Steiner. As an approach to this we prove the following. 

\begin{theorem} \label{theo:zigzag} Every graph $G$ with
    $\zig(G)=\chi(G)= t$ contains $K_{\piso{\frac t2}+1}$ as a totally odd immersion.
\end{theorem}

This result, along with Theorem~\ref{theo:Kneser2}, give evidence toward the following conjecture $-$ given in~\cite{JimenezQuirozThraves} and inspired by a weaker conjecture of Churchley~\cite{C17} $-$  which can be seen as the immersion-analogue of the Odd Hadwiger Conjecture.

\begin{conjecture}[Jim\'enez, Quiroz, Thraves Caro]\label{conj:ours}
    Every graph $G$ contains $K_{\chi(G)}$ as a totally odd immersion.
\end{conjecture}

Conjecture~\ref{conj:ours} holds trivially when $\chi(G)\le 3$ and, by a result of Thomassen~\cite{Thomassen} and (independently) Zang~\cite{Z98}, it also holds for $\chi(G)=4$. In light of Theorem~\ref{theo:zigzag} we make the following weakening of Conjecture~\ref{conj:ours}.

\begin{conjecture}\label{conj:new}
    There exist $\epsilon > 0$ such that every graph $G$ contains $K_{\lceil \epsilon \chi(G)\rceil}$ as a totally odd immersion.
\end{conjecture}

After the appearance of a first draft of this paper, McFarland~\cite{McFarland} proved an analogous result for a weak notion of immersions. But still the best bound known in general for totally odd (strong) immersions, is implied by a result of Kawarabayashi~\cite{K13}, who proved that every graph $G$ contains a totally odd subdivision of $K_{\sqrt{\frac 4{79}\chi(G)}}$.

Our last result is concerned with totally odd substructures in Generalised Mycielski graphs\footnote{Given a graph $G$ and a integer $m\geq 2$, the \emph{$m$-level generalised Mycielskian} of $G$, denoted $\mu_m(G)$, is the graph with vertex set $(V(G) \times \{0,1,\dots,m-1\}) \cup \{w\}$, and edges $(u,0)(v,0)$ and $(u,i)(v,i+1)$ for every $uv \in E(G)$ and $i\in\{ 0,1,\dots ,m-2\}$, together with edges $(u,m-1)w$ for all $u\in V(G)$.}, which are another classic example of graphs having chromatic number equal to its topological bounds. It was Gy\'arf\'as, Jensen and Stiebitz~\cite{GyarfasJenseStiebitz} who showed that these graphs have chromatic number equal to its topological bounds, while Tardif~\cite{ClaudeTardif} showed that the lower bound $\lceil \frac{|V(G)|}{\alpha(G)} \rceil$ can be far from the chromatic number, for such a graph $G$. 

It is not hard to see that $\chi(\mu_m(G))\le \chi(G)+1$ for every $m\ge 2$. Mycielski~\cite{Mycielski} observed that $\chi(\mu_2(G))=\chi(G)+1$ for every graph $G$, while Tardif~\cite{ClaudeTardif} showed that for every $t\ge 4$ there exists $H$ with $t=\chi(H)=\chi(\mu_3(H))$. We show the following, which discards Generalised Mycielski graphs as minimal counterexamples to Conjecture~\ref{conj:ours}.

\begin{theorem}\label{theo:Mycielski}
    Let $m \geq 2$ be an integer. If $G$ contains $K_t$ as a totally odd immersion, then $\mu_m(G)$  contains $K_{t+1}$ as a totally odd immersion. 
\end{theorem}

Our proof technique also yields the following analogous result for totally odd subdivisions.

\begin{theorem}\label{theo:Mycielski2}
    Let $m \geq 2$ be an integer. If $G$ contains $K_t$ as a totally odd subdivision, then $\mu_m(G)$  contains $K_{t+1}$ as a totally odd subdivision. 
\end{theorem}

Analogous results have been obtained by Collins, Heenehan, and McDonald~\cite{CollinsHeenehanMcDonald} and  Simonyi and Zsb\'an \cite{SimonyiZsban} for immersions and odd minors, respectively.

An extended abstract containing sketches of the proofs of Theorem~\ref{theo:Kneser} and Theorem~\ref{theo:zigzag} is to appear in \cite{us}.

The rest of the paper is organised as follows. In Section \ref{sec:arbitraryKneser} we prove Theorem \ref{theo:Kneser}, and in Section \ref{sec:chikneser} we prove Theorem~\ref{theo:Kneser2}. In Section~\ref{sec:toizig} we prove Theorem~\ref{theo:zigzag}, and in Section~\ref{sec:submicyelski} we prove Theorem~\ref{theo:Mycielski}. 

We end the introduction with some notation. For any integer $n\ge1$, we set $[n]:=\{1,2,\dots ,n\}$. If $a\leq b$ are integers, we set $[a,b] :=\{a,a+1,\ldots,b-1,b\}$. If $b<a$, then $[a,b]$ is the empty set.

\section{Proof of Theorem \ref{theo:Kneser}}  \label{sec:arbitraryKneser} 

\noindent
Let $k\ge 1$ and $r\ge 6$ be integers. For $A\subset [2k+r]$ we say that $a\in A$ is an \emph{isle} of $A$ if $\{a-1,a,a+1\}\cap A=\{a\}$, and that $a$ is a \emph{gap} of $A$ if $\{a-1,a,a+1\}\cap A=\{a-1,a+1\}$. We set $$\mathcal{I}_{k,r}= \left\{ i \in \mathbb{N}\mid \frac{k-1}{2} + 1 \leq i\leq k+r{-3} \,\,\,\textrm{ and } \,\, i\equiv0 \mod 4 \right\},$$ and note that our choice of $r$ gives $|\mathcal{I}_{k,r}|\ge (k+r)/8$. It is well known and not hard to see that the chromatic number of the Kneser graph $KG(2k+r,k)$ is at most $r+2$ (which Lov\'asz \cite{Lovasz} showed to be a tight bound). Therefore, to prove Theorem~\ref{theo:Kneser} it suffices to show the following.

\begin{theorem}\label{theo:kover8}
    Let $k\ge 13$ and $r\ge 6$ be integers, with $k-1$ divisible by $2r$. Then the Kneser graph $KG(2k+r,k)$ contains $K_{\lceil \frac{k+r}{8} \rceil}$ as a totally odd subdivision.
\end{theorem}

Before we proceed with the proof, we mention that we believe Theorem~\ref{theo:kover8} can be improved to obtain larger complete graphs as totally odd subdivisions, and to weaken the divisibility restriction, but we do not attempt these improvements for the sake of brevity and readability.

\begin{proof}
    As it suffices, we will construct a totally odd subdivision of a complete graph on $|\mathcal{I}_{k,r}|$ vertices. As terminals we choose the vertices
$$X_i:=\{1,\ldots,k-1\}\cup\{k+i\}$$
 with $i \in \mathcal{I}_{k,r}$. To find the path joining terminals $X_i$ and $X_j$ we are going to define
two sets of vertices,
$B_{i,j}^1,\ldots,B_{i,j}^{(k-1)/r}$ and
$C_{i,j}^1,\ldots,C_{i,j}^{(k-1)/r},$
such that
$$P_{ij}:=X_i,B_{ij}^1,\ldots ,B_{ij}^{(k-1)/r},C_{ij}^{(k-1)/r},\ldots ,C_{i,j}^1,X_j$$
is a path. Notice that the length of such a path would be 
$2\frac{k-1}{r} -1$. The definitions of $B_{ij}^s$ and $C_{ij}^s$ will depend on the parity of $s$ and, in particular, this will help guarantee the adjacency between $B_{ij}^s$ and $B_{ij}^{s+1}$ ($C_{ij}^s$ and $C_{ij}^{s+1}$ resp.).

For any fixed pair $i,j \in \mathcal{I}_{k,r}$ with $i<j$, our main goal is to define the vertices of the $(X_i,X_j)$-path of our subdivision. To distinguish the $(X_i,X_j)$-path from other paths, the vertices of this path will contain $k+i-1$ and $k+j$ as either their last two isles or their last two gaps. So first we set $S_{i,j}=\{k+i-1,k+j\}$ and $T_{i,j}=\{k+i-2,k+i,k+j-1,k+j+1\}.$ Whenever $p$ is such that $2p+1\le (k-1)/r$ we define
$$
B_{i,j}^{2p+1}=[1,pr]\cup S_{i,j}\cup [k+(p+1)r-3,k+i-3] \cup [k+i+1,k+j-2]  \cup [k+j+2,2k+r],
$$
and whenever $2p+2\le (k-1)/r$ we define
$$
B_{i,j}^{2p+2}=[(p+1)r+1,k+(p+1)r-4]\cup T_{i,j}.
$$

Let us now study the isles and gaps of $B_{i,j}^{2p+1}$ and $B_{i,j}^{2p+2}$. Since we have $i,j \in \mathcal{I}_{k,r}$ and $\frac{k-1}{2r}\in \mathbb{N}$, for all possible values of $p$, each of the intervals $[k+(p+1)r-3,k+i-3]$ , $[k+i+1,k+j-2]$ and $  [k+j+2,2k+r]$ has at least two elements. Furthermore, either $p=0$ and hence $[1,pr]$ is empty, or $p\geq 1$ and hence $|[1,pr]|\geq r\geq 2$. So the only isles of  $B_{i,j}^{2p+1}$ are $k+i-1$ and $k+j$. Now for $B_{i,j}^{2p+2}$, we first see that the inequality $k+(p+1)r-2 < k+i-2$ tells us that $k+(p+1)r-3$ is not a gap, while the fact that $j-i\geq 4$ give us that there is no gap between $k+i$ and $k+j-1$. Thus the only gaps in $B_{i,j}^{2p+2}$ are $k+i-1$ and $k+j$. Also by definition, successive sets $B_{i,j}^{2p+1}$, $B_{i,j}^{2p+2}$, and  $B_{i,j}^{2p+2}$, $B_{i,j}^{2p+3}$ are disjoint for all relevant values of $p$. 

Now we define the sets ${C}_{i,j}^{2p+1}$ and ${C}_{i,j}^{2p+2}$ depending on the value of $p$.

\noindent
$\bullet$ If $pr\geq k-i-1$, we define
$$C^{2p+1}_{i,j}=[k-pr,2k-pr-5]\cup T_{i,j}$$
and
 $$C^{2p+2}_{i,j}=[1,k-(p+1)r-1]\cup S_{i,j}\cup [2k-pr-4,k+i-3]\cup [k+i+1,k+j-2]\cup [k+j+2,2k+r].$$
According to our definition, $C^{2p+1}_{i,j}$ contains $k+i-1$ and $k+j$ as the last two gaps since we have $2k-pr-3 \leq k+i-2$. Note that in the case that $pr= k-i-1$, $k+i-3$ is also a gap (though smaller than $k+i-1$ and $k+j$).
Now, let us check that  $k+i-1$ and $k+j$ are the last two isles of $C^{2p+2}_{i,j}$. Note that in the case $pr= k-i-1$, the element $2k-pr-4=k+i-3$ is an isle, though smaller than $k+i-1$, and for $pr\geq k-i$, the interval $[2k-pr-4,k+i-3]$ has at least two elements. Furthermore, the intervals  $[1,k-(p+1)r-1]$, $[k+i+1,k+j-2]$  and $[k+j+2,2k+r]$ all have at least two elements (respectively, because we always have $2p+1\le (k-1)/r$, we also always have $j-i \geq 4$, and by definition of $\mathcal{I}_{k,r}$). 

Now, if $pr\leq k-i-2$, we split the definition into the following three other cases depending on the relation between $p$ and $j$. The analysis in each case is analogous to the previous ones.

\noindent
$\bullet$ If  $pr\leq k-j$, we define
$$C^{2p+1}_{i,j}=[k-pr,k+i-2]\cup [k+i, k+j-1] \cup [k+j+1, 2k-pr+1]$$
and
 $$C^{2p+2}_{i,j}=[1,k-(p+1)r-1]\cup S_{i,j} \cup[2k-pr+2, 2k+r].$$
In this case, the only two gaps in $C^{2p+1}_{i,j}$ are $k+i-1$ and $k+j$ since each of the three intervals forming the set has at least two elements. For  $C^{2p+2}_{i,j}$,  $r\geq 3$ gives us that the interval $[2k-pr+2, 2k+r]$ has at least two elements for all relevant values of $p$. The interval $[1,k-(p+1)r-1]$ also has at least two elements. We have $k+j \leq 2k-pr$ because of the case assumption, and $k-(p+1)r+1 \leq k+i-3$. Altogether, we have that the only two isles of $C^{2p+2}_{i,j}$ are $k+i-1$ and $k+j$.

\noindent
$\bullet$  If  $pr\geq k-j+2$ we define
$$C^{2p+1}_{i,j}=[k-pr,k+i-2]\cup [k+i,2k-pr-2]\cup \{k+j-1,k+j+1\}$$
and
 $$C^{2p+2}_{i,j}=[1,k-(p+1)r-1]\cup S_{i,j}\cup [2k-pr-1,k+j-2]\cup [k+j+2,2k+r].$$
In this case, $k+i-1$ and $k+j$ are the only two gaps of $C^{2p+1}_{i,j}$, and the only two isles of $C^{2p+2}_{i,j}$.

\noindent
$\bullet$ If $pr= k-j+1$, we define
$$C^{2p+1}_{i,j}=[k-pr,k+i-2]\cup [k+i,2k-pr-3]\cup \{k+j-1,k+j+1,k+j+2\}$$
and
 $$C^{2p+2}_{i,j}=[1,k-(p+1)r-1]\cup S_{i,j}\cup [2k-pr-2,k+j-2]\cup [k+j+3,2k+r].$$
Note that since $j-i \geq 4$, we have  $pr\leq k-i-3$. 
This helps us see that $k+i-1$ and $k+j$ are the only two gaps of $C^{2p+1}_{i,j}$, while the only two isles of $C^{2p+2}_{i,j}$.

It is easy to check that, with the given definitions, successive sets $C_{i,j}^{2p+1}$ and $C_{i,j}^{2p+2}$ are disjoint for all relevant values of $p\geq 0$. In order to show that $C_{i,j}^{2p+2}$ and $C_{i,j}^{2p+3}$ are disjoint, we note that  
$${C}_{i,j}^{2p+2} \subseteq ([1,k-(p+1)r-1]\cup [2k-pr-4,2k+r]\cup S_{i,j})\setminus T_{i,j},$$
in each case,  while 
\begin{align*}
{C}_{i,j}^{2p+3}={C}_{i,j}^{2(p+1)+1} \subseteq &([k-(p+1)r,2k-(p+1)r+1]\cup T_{i,j})\setminus S_{i,j}\\ 
= & ([k-(p+1)r,2k-pr-(r-1)]\cup T_{i,j})\setminus S_{i,j}
\end{align*}
in every case except for that in which $(p+1)r=k-j+1$, where $k+j+2$ is also contained in ${C}_{i,j}^{2p+3}$.
Thus, whenever  $(p+1)r\neq k-j+1$, we have that $C_{i,j}^{2p+2}$ and $C_{i,j}^{2p+3}$ are disjoint. For the case $(p+1)r= k-j+1$, we have
\begin{align*}
C^{2(p+1)+1}_{i,j}=&[k-(p+1)r,k+i-2]\cup [k+i,2k-(p+1)r-3]\cup \{k+j-1,k+j+1,k+j+2\}
\\
=&[k-(p+1)r,k+i-2]\cup [k+i,k+j-4]\cup \{k+j-1,k+j+1,k+j+2\}
\end{align*}
and (taking $pr=k-j+1-r$) 
 \begin{align*} 
 C^{2p+2}_{i,j}=&[1,k-(p+1)r-1]\cup S_{i,j} \cup[2k-pr+2, 2k+r]\\
 =& [1,k-(p+1)r-1]\cup S_{i,j} \cup[k+j+r+1, 2k+r]
\end{align*}
and hence, $C_{i,j}^{2p+2}$ and $C_{i,j}^{2p+3}$ are disjoint as well in this case.

To finish checking that consecutive vertices of $P_{ij}$ are indeed adjacent, we note that when we have $2p+2=\frac{k-1}{r}$, the following inequalities hold
$$k+(p+1)r-4<2k-pr-4
\,\, \text{and} \,\, 
 k-(p+1)r-1<(p+1)r+1,$$
which gives us 
$${B}_{ij}^{(k-1)/r} \cap {C}_{ij}^{(k-1)/r} =\varnothing.$$

Now to see that $P_{ij}$ is indeed a path, note that by construction, $B_{ij}^1,\ldots ,B_{ij}^{(k-1)/r}$ are all distinct; the same for $C_{ij}^1,\ldots ,C_{i,j}^{(k-1)/r}$. Suppose then that there are values $s$, $s'$ such that $B_{ij}^{s} = C_{ij}^{s'}$. 
For such an equality, we have two possibilities depending on whether $1\in B_{ij}^{s}= C_{ij}^{s'}$ or not. If $1$ is in said sets, then we must have $s=2p+1$ and $s'=2p'+2$ for some $0\leq p,p'\leq \frac{k-1}{2r}-1$. But the equality of the sets implies $[1,pr]=[1,k-(p'+1)r-1]$, which implies $k=(p+p'+1)r+1$, contradicting $p,p'\leq \frac{k-1}{2r}-1$. On the other hand, if $1\not\in B_{ij}^{s}= C_{ij}^{s'}$, we must have $s=2p+2$ and $s'=2p'+1$. In such a case, the smallest element in $B_{ij}^{s}$ is $(p+1)r+1$ and the smallest element in $C_{ij}^{s'}$ is $k-p'r$. So $B_{ij}^{s}= C_{ij}^{s'}$ implies $k=(p+p'+1)r+1$, contradicting again $p,p'\leq \frac{k-1}{2r}-1$. 

To finish the proof, let us see that the paths are mutually disjoint. Consider two different pairs of vertices $X_i, X_j$ and $X_h, X_\ell$. Suppose for a contradiction that there exists a vertex $Z$ at the intersection of $P_{ij}$ with $P_{h\ell}$.  Assume first that $Z\cap [1,k-1]\ne \varnothing$. If we further have $1\in Z$, then we must have $Z=B_{ij}^{2p+1}$ or $Z=C_{ij}^{2p+2}$ for some value of $p$. Either way, the two last isles of~$Z$ force $\{i,j\}=\{h,\ell\}$, a contradiction. If we instead have $1\not\in Z$, then it must be $Z=B_{ij}^{2p+2}$ or $Z=C_{ij}^{2p+1}$ for some value of $p$. In this case, their two last gaps force $\{i,j\}=\{h,\ell\}$. Then we can assume that $Z\cap [1,k-1]=\varnothing$, in which case we must have $Z=B_{ij}^1$ or  $Z=C_{ij}^1$. If $|Z\cap [2k+2,2k+r]|\geq 5$, then $Z=B^1_{ij}$, and we can again look at the isles to obtain $\{i,j\}=\{h,\ell\}$. In the case that $|Z\cap [2k+2,2k+r]|\leq 4$, we must have $Z=C^1_{ij}$, and we can look at the gaps to get $\{i,j\}=\{h,\ell\}$, and thus a contradiction. It follows that the paths joining the terminals are mutually vertex-disjoint. The result follows.
\end{proof}

\section{Proof of Theorem \ref{theo:Kneser2}}\label{sec:chikneser} 

As mentioned in the previous section,
the chromatic number of the Kneser graph $\Gamma=KG(2k+r,k)$ is $\chi(\Gamma)=r+2$, and thus we want to find a totally odd subdivision of $K_{r+2}$. Since Thomassen~\cite{Thomassen} and Zang~\cite{Z98} proved that every 4-chromatic graph contains $K_4$ as a totally odd subdivision, we can restrict ourselves to $r\geq 3$ (and since the result is otherwise trivial, we take $k\ge 2$). 

\subsection{When $r\geq k$}
In this case we actually construct a subdivision where each path joining two terminals has length either 1 or 3. We present the subcases $k=2,3$ in the appendix, and work here with $k\geq 4$.
Recall that we are looking for a totally odd  subdivision of $K_{r+2}$. 
We chose as terminals the vertices
\begin{align*}
Y=&\{1,\ldots,k\},\\
X_i=&\{k+1,\ldots,2k-1\}\cup\{2k+i\},
\end{align*}
with $0\leq i\leq r$.

For fixed $0<i<j<r$, we join $X_i$ to $X_j$, through a path $X_iB_{i,j}C_{i,j}X_j$, whose internal vertices are defined according to three cases which depend on the relationship between $i,j$, and $k$. If $j-i\ge k-2$, we choose $a,b$ such that $j-i-(b-a)=k-3$ and
\begin{align*}
B_{i,j}=&[1,k-2]\cup\{2k+i-1,2k+j\},\\
C_{i,j}=&[k-1,k]\cup[2k+i,a]\cup[b,2k+j-1],
\end{align*}
are sets of size $k$. Otherwise, if $i\ge 2$ and $j-i\le k-4$, choose $a,b$ such that $b-a=k-1$, $2k\le a\le 2k+i-2$, $2k+j+1\le b\le 2k+r$ and set
\begin{align*}
B_{i,j}=&[1,k-2]\cup\{2k+i-1,2k+j\},\\
C_{i,j}=&[k-1,k]\cup[a,2k+i-2]\cup[2k+i,2k+j-1]\cup[2k+j{+1},b].
\end{align*} Finally, if $i=1$ or $j-i=k-3$, we define \begin{align*}
    B_{i,j}=&[1,k-2]\cup\{2k+i-1,2k+j\},\\
    C_{i,j}=&[k-1,k]\cup [2k+i,2k+j-1]\cup [2k+j+1,b],
\end{align*}
where $b\leq 2k+r$ is such that $|C_{i,j}|=k$. 

We now define the internal vertices of the $(X_0,X_j)$-paths for $j\in\{1,\dots r-1\}$.  If $j<k$, set
\begin{align*}
B_{0,j}=&[1,k-1]\cup\{2k+j\}\\
C_{0,j}=&\{k\}\cup[2k+1,2k+j-1]\cup[2k+j+1,3k].
\end{align*}
Otherwise, let $2k\leq a<2k+j<b\leq 2k+r$ be such that $b-a=k-1$, and set 
\begin{align*}
B_{0,j}=&[1,k-1]\cup\{2k+j\}\\
C_{0,j}=&\{k\}\cup[a,2k+j-1]\cup[2k+j+1,b].
\end{align*}

To finish, we need to give the paths to $X_r$. If $r-k+2\le i\le r-1$, set 
\begin{align*}
B_{i,r}=&[1,k-2]\cup\{2k+i-1,2k+r\}\\
C_{i,r}=&[k-1,k]\cup[k+r+1,2k+i-2]\cup[2k+i,2k+r-1].
\end{align*}
If $i<r-k+2$, set 
\begin{align*}
B_{i,r}=&[1,k-3]\cup\{2k+i-1,2k+r-1,2k+r\}\\
C_{i,r}=&[k-1,k]\cup\{2k+i-2,2k+i\}\cup[k+r+3,2k+r-2].
\end{align*}
Moreover, set 
\begin{align*}
B_{0,r}=&[1,k-1]\cup\{2k+r\}\\
C_{0,r}=&\{k\}\cup[k+r+1,2k+r-1].
\end{align*}

Now notice that, in any case, $X_iB_{i,j}C_{i,j}X_j$ is indeed a path, and we must check that different paths are internally disjoint. Indeed, $1\in B_{ij}$ for every $\{i,j\}$ but $1\not\in C_{h,\ell}$ for every $\{h,\ell\}$, and thus $B_{ij}\neq C_{h\ell}$ always holds. Moreover, vertices $B_{ij}$ are determined by the elements $2k+i-1$ and $2k+j$, while vertices $C_{ij}$ are determined by the elements $2k+i$ and $2k+j-1$, hence different pairs $\{i,j\}\neq \{h,\ell\}$ determine different sets $B_{ij}\neq B_{h,\ell}$ and $C_{ij}\neq C_{h,\ell}$. Thus, we have the desired subdivision.

\subsection{When $r< k$}
Consider the graph $KG(2k+r,k)$ for some fixed choice of $r$ and $k$ such that $3\le r< k$. We  are looking for a totally odd subdivision of $K_{r+2}$, and choose for terminals the vertices
\begin{align*}
Y=&\{1,\ldots,k\}\\
X_i=&\{k+1,\ldots,2k-1\}\cup\{2k+i\},
\end{align*}
for $0\leq i\leq r$. In a manner similar to the proof of Theorem~\ref{theo:kover8}, every $(X_i, X_j)$-path will have length
$2\lceil(k-1)/r\rceil +1$ and be of the form
\[
X_i,B_{i,j}^1,\ldots ,B_{i,j}^{\lceil(k-1)/r\rceil},C_{i,j}^{\lceil(k-1)/r\rceil},\ldots ,C_{i,j}^1,X_j
\]
for some adequate choice of inner vertices. As $B_{i,j}^s$ must be disjoint from $B_{i,j}^{s+1}$, their definition will depend on the parity of $s$; the same will happen with $C_{i,j}^s$. 
Although we give a general definition for $B_{i,j}^s$ and $C_{i,j}^s$, we only consider the cases when $1\leq s\leq \lceil (k-1)/r\rceil$. We will also
split the definition in two cases, whether $i=0$, or $1\le i<j\leq r$. 

Assume $1\le i<j\leq r$ and for $p\geq 0$ let
\begin{align*}
B_{i,j}^{2p+1}=&[1,k-2-pr]\cup[k+1,k+pr]\cup\{2k+i-1,2k+j\}\\
B_{i,j}^{2p+2}=&[k-1-pr,k]\cup[k+(p+1)r+1,2k+i-2]\cup[2k+i,2k+j-1]\cup[2k+j+1,2k+r]\\
C_{i,j}^{2p+1}=&[(p+1)r,k]\cup[2k-pr,2k+i-2]\cup[2k+i,2k+j-1]\cup[2k+j+1,2k+r]\\
C_{i,j}^{2p+2}=&[1,(p+1)r-1]\cup[k+1,2k-(p+1)r-1]\cup\{2k+i-1,2k+j\},
\end{align*}
Notice that both $2p+1\leq \lceil (k-1)/r\rceil$, $2p+2\leq \lceil (k-1)/r\rceil$, imply $pr\leq (k-1)/2$.
Thus is not hard to see that for every relevant value of $s$  we have $B_{i,j}^s\cap B_{i,j}^{s+1}=C_{i,j}^s\cap C_{i,j}^{s+1}=\emptyset$, giving us that $B_{i,j}^1,\ldots ,B_{i,j}^{\lceil(k-1)/r\rceil}$ and $C_{i,j}^1,\ldots ,C_{i,j}^{\lceil(k-1)/r\rceil}$ are paths. 
We need to show now that $B_{i,j}^{\lceil(k-1)/r\rceil}$ is disjoint from $C_{i,j}^{\lceil(k-1)/r\rceil}$.
To do this, notice that
\begin{align*}
B_{i,j}^{2p+1}\cap C_{i,j}^{2p+1}=&\begin{cases}
[(p+1)r,k-2-pr]&\text{if $(2p+1)r\leq k-2$}\\
\emptyset &\text{if $2pr< k-1\leq (2p+1)r$}\end{cases}\\
B_{i,j}^{2p+2}\cap C_{i,j}^{2p+2}=&\begin{cases}[k+(p+1)r+1,2k-(p+1)r-1]&\text{if $(2p+2)r\leq k-2$}\\
\emptyset &\text{if $(2p+1)r< k-1\leq (2p+2)r$.}\end{cases}
\end{align*}

In the first case, the intersection is empty when $(p+1)r>k-2-pr$, that is, when $pr>\frac{k-r-2}{2}$. This is the case when $2p+1=\lceil\frac{k-1}{r}\rceil$. Similarly, for the even case, the intersection is empty when $k+(p+1)r+1>2k-(p+1)r-1$, i.e., when $pr>\frac{k-2-2r}{2}$. Which again holds particularly when $2p+2=\lceil\frac{k-1}{r}\rceil$. Then,
$X_i,B_{i,j}^1\ldots , B_{i,j}^{\lceil(k-1)/r\rceil},C_{i,j}^{\lceil(k-1)/r\rceil},\ldots ,C_{i,j}^1,X_j$
 is a path.

We need now to define the paths from $X_0$ to all other $X_j$'s. To do this, we let
\begin{align*}
B_{0,j}^{2p+1}=&[1,k-1-pr]\cup[k+1,k+pr]\cup\{2k+j\}\\
B_{0,j}^{2p+2}=&[k-pr,k]\cup[k+(p+1)r+1,2k+j-1]\cup[2k+j+1,2k+r]\\
C_{0,j}^{2p+1}=&[(p+1)r+1,k]\cup[2k-pr,2k+j-1]\cup[2k+j+1,2k+r]\\
C_{0,j}^{2p+2}=&[1,(p+1)r]\cup[k+1,2k-(p+1)r-1]\cup\{2k+j\}.
\end{align*}
As before, we get $B_{0,j}^{2p+1}\cap B_{0,j}^{2p+2}=B_{0,j}^{2p+2}\cap B_{0,j}^{2p+3}=\emptyset$ and
$C_{0,j}^{2p+1}\cap C_{0,j}^{2p+2}=C_{0,j}^{2p+2}\cap C_{0,j}^{2p+3}=\emptyset$ for all relevant values of $p$. We also get
\begin{align*}
B_{0,j}^{2p+1}\cap C_{0,j}^{2p+1}=&\begin{cases}
[(p+1)r+1,k-1-pr]&\text{if $(2p+1)r\leq k-2$}\\
\emptyset &\text{if $2pr< k-1\leq (2p+1)r$}\end{cases}\\
B_{0,j}^{2p+2}\cap C_{0,j}^{2p+2}=&\begin{cases}[k+(p+1)r+1,2k-(p+1)r-1]&\text{if $(2p+2)r\leq k-2$}\\
\emptyset &\text{if $(2p+1)r< k-1\leq (2p+2)r$,}\end{cases}
\end{align*}
which give us a path of length $2\lceil (k-1)/r\rceil +1$ from $X_0$ to $X_j$.

We need to show that the vertices in the $(X_i,X_j)$-path are different from the vertices in the $(X_h,X_\ell)$-path whenever $\{i,j\}\neq \{h,\ell\}$. To do this, let $Z\in\{B_{i,j}^s,C_{i,j}^s\}$, for some $0\leq i<j\leq r$, $1\leq s\leq \lceil(k-1)/r\rceil$. Consider whether $1\in Z$ or $1\not\in Z$. If $1\in Z$, then $Z$ can only be $B_{i,j}^{2p+1}$ or $C_{i,j}^{2p+2}$ for some $i,j,p$. We can determine $i,j$ by considering $Z\cap [2k+1,2k+r]$. Notice that $Z\cap [2k+1,2k+r]$ is $\{2k+j\}$  if and only if $i=0$, and it is 
$\{2k+i-1,2k+j\}$  otherwise. 
We can now determine whether $Z=B_{i,j}^{2p+1}$ or $Z=C_{i,j}^{2p+2}$ by  considering $[1,k]\cap Z=[1,b]$. Notice that the possible values of $b$ are quite restricted, i.e., $b\in \{k-2-pr,k-1-pr\}$ if $Z=B_{i,j}^{2p+1}$ and $b\in\{(p'+1)r-1,(p'+1)r\}$ if $Z=C_{i,j}^{2p'+2}$. Thus, we have
\[
Z=\begin{cases}
B_{i,j}^{2[(k-2-b)/r]+1}&\text{if $2[(k-2-b)/r]+1\leq \lceil(k-1)/r\rceil$ and $i\neq 0$,}\\
B_{0,j}^{2[(k-1-b)/r]+1}&\text{if $2[(k-1-b)/r]+1\leq \lceil(k-1)/r\rceil$ and $i=0$,}\\
C_{i,j}^{2[-1+(b+1)/r]+2}&\text{if $2[-1+(b+1)/r]+2\leq \lceil(k-1)/r\rceil$ and $i\neq 0$},\\
C_{0,j}^{2[-1+(b/r)]+2}&\text{if $2[-1+(b/r)]+2\leq \lceil(k-1)/r\rceil$ and $i=0$.}\\
\end{cases}
\]
Let us see that it can be only one of those cases. 
Assume otherwise. As we already know the value of $i$, we have two options. Either $i\neq 0$, $b=k-2-pr$ for some $p$ and $b=(p'+1)r-1$ for some $p'$, or $i=0$, $b=k-1-pr$ for some $p$ and $b=(p'+1)r$  for some $p'$. If $i\neq 0$  notice that $b=k-2-pr$ implies that
$r\mid (k-2-b)$ whereas $b=(p'+1)r-1$ implies $r|(b+1)$. Thus $r$ divides $k-1$, which implies $\lceil (k-1)/r\rceil=(k-1)/r$. But 
\[
2[(k-2-b)/r]+1+2[-1+(b+1)/r]+2=1+2(k-1)/r=1+2\lceil(k-1)/r\rceil,
\]
contradicting the fact that $2p+1,2p'+2\leq \lceil(k-1)/r\rceil$. 

On the other hand, when $i=0$ $b=k-1-pr$  implies that $r\mid k-1-b$ whereas $b=(p'+1)r$ implies $r\mid b$, which implies that $r|(k-1)$. But
\[
2[(k-1-b)/r]+1+2[-1+(b/r)]+2=1+2(k-1)/r=1+2\lceil(k-1)/r\rceil,
\]
contradicting again the upper bounds. 
Thus, if $1\in Z$ we can determine whether it is a $B^{2p+1}_{i,j}$ or a $C_{i,j}^{2p+2}$ set, and what the value of $i,j,$ and $p$ are.

If $1\not\in Z$, then $Z$ is either $B_{i,j}^{2p+2}$ or $C_{i,j}^{2p+1}$. To decide which, consider $Z\cap [k+1,2k]=[b,2k]$. Again, notice that $b=k+(p+1)r+1$ if $Z=B_{i,j}^{2p+2}$ and $b=2k-pr$ if $Z=C_{i,j}^{2p+1}$. Thus, we have
\[
Z=\begin{cases}
B_{i,j}^{2[-1+(b-k-1)/r]+2}&\text{if $2[-1+(b-k-1)/r]+2\leq \lceil(k-1)/r\rceil$,}\\
C_{i,j}^{2[(2k-b)/r]+1}&\text{if $2[(2k-b)/r]+1\leq \lceil(k-1)/r\rceil$}.\\
\end{cases}
\]
Notice that it can only be one of those cases, because if both held 
we would have $r\mid b-k-1$ and $r\mid (2k-b)$, which implies $r\mid k-1$ and gives
\[
2[-1+(b-k-1)/r]+2+2[(2k-b)/r]+1=2[(k-1)/r]+1=2\lceil (k-1)/r\rceil+1,
\]
contradicting the upper bounds.
Thus, if $1\not\in Z$ we can determine whether it is a $B$ or a $C$ set, and what the value of $p$ is. Further, in this case, let $\overline{Z}$ be the complement of $Z$, and notice that  
$\overline{Z}\cap [2k+1,2k+r]$ is $\{2k+j\}$  if and only if $i=0$, and it is 
$\{2k+i-1,2k+j\}$  otherwise. Thus, we can recognize $i$ and $j$ from $\overline{Z}\cap [2k+1,2k+r]$. 

As $Z$ can be completely determined, the paths obtained are pairwise internally vertex-disjoint.

\section{Odd immersions in graphs with $\zig(G)=\chi(G)$}\label{sec:toizig}

\noindent
In this section, we prove Theorem~\ref{theo:zigzag} which says that a graph with both zigzag and chromatic number equal to $t$ contains $K_{\piso{\frac t2}+1}$ as a totally odd immersion.  

Our approach here is to take a maximum zigzag of $(G,c)$, for some optimal colouring $c$ of $G$, and to consider one of the sides of this complete bipartite graph as the terminals of our totally odd immersion. The edges of this bipartite graph would immediately give us a complete immersion, but it would be totally even. So we would instead try to connect our terminals through Kempe chains (components of the graph induced in $G$ by those vertices coloured $i$ or $j$, for some fixed choice of $i,j$). Kempe chains are edge disjoint when different pairs of colours are considered, a fact which has already been used in the construction of immersions in \cite{Biclique} and~\cite{GLW19}. Moreover, if a path in a Kempe chain starts and finishes in vertices of different colours, then this path is odd. So whenever two terminals are in the same Kempe chain, we have the desired odd path joining them. If a terminal $v_i$ coloured $i$ is not in the same $\{i,j\}$-Kempe chain as the terminal $v_j$ coloured $j$, then (with the additional assumption that $c$ is taken so as to minimize the number of maximum zigzags) $v_i$ is in the same component as a vertex which lies in a complete bipartite graph which is one Kempe switch away from becoming a maximum zigzag. From this ``potential zigzag'' departs some other $\{i,j\}$-Kempe chain, which we can follow to another such complete bipartite graph. Our first goal is to show that, at some point, this process leads us back to our first maximum zigzag. This gives us a series of paths between zigzags and potential zigzags, which starts and ends at the vertices we want to join. Then we want to sew together this collection of paths into one (odd) path joining our desired pair of vertices. To do this we need edges that are not in the $\{i,j\}$-Kempe chain, and so we need to be careful in our choice of this edges so that, in the end, we can guarantee that all paths are edge disjoint. 

We start with an auxiliary lemma.

\begin{lemma}\label{lemma:existevennumber}
Let $t\ge5$ be an integer. For every pair of distinct odd $i,j\in [t]$, we can choose two distinct even numbers $\ell_{\{i,j\}}^i, \ell_{\{i,j\}}^j \in [t-1]$ such that for any odd $k \in[t]\setminus \{i,j\},$ we have $\ell_{\{i,j\}}^i \neq \ell_{\{i,k\}}^i$.
\end{lemma}
\begin{proof}
    For every pair of distinct odd $i,j\in [t]$, we define dummy variables $x_{\{i,j\}}^i$ and $x_{\{i,j\}}^j$ such that all the variables are distinct. Let $H$ be the graph having
$\conjunto{x_{\{i,j\}}^i \mid i, j\in[t] \text{ odd and distinct}}$ as its vertex set and $\{x_{\{i,j\}}^i x_{\{i,j\}}^j, x_{\{i,j\}}^i x_{\{i,k\}}^i \mid i,j,k \in [t] \text{ are distinct odd numbers}\}$ as its edge set. 
It is not hard to see that $H$ is connected and $(\lceil t/2 \rceil -1)$-regular (the neighbourhood of $x_{\{i,j\}}^i$ is exactly $x_{\{i,j\}}^j$ together with each $x_{\{i,k\}}^i$ for $k \in [t]$ odd and different from $i$ and $j$). Thus by Brook's Theorem, there exists a $(\lceil t/2 \rceil -1)$-colouring of $H$. Let $c$ be such a colouring with its colours being the even numbers of $[t-1]$. Now define $\ell_{i,j}^i := c(x_{\{i,j\}}^i)$, for every pair of odd $i,j \in [t]$. The result follows.
\end{proof}

We now obtain a couple of observations and a lemma for graphs $G$ with $\zig(G)=\chi(G)=t$. For these we need some definitions and notations.

For a graph $G$ and a proper colouring $c$ of $G$, 
we use $G_{ij}$ to denote the graph induced in $G$ by vertices coloured $i$ or $j$ by $c$. For a vertex $x\in V(G)$ with $c(x)\in\{i,j\}$, we let $C_{ij}^x$ be the component of $G_{ij}$ that contains $x$ (that is, the $\{i,j\}$-Kempe chain containing $x$).
Let $Z=z_1,z_2, \dots ,z_t$ be a sequence of (distinct) vertices of $G$. We say that $Z$ is an \emph{$\{i,j\}$-potential zigzag} in $(G,c)$ if there is a pair of vertices $u,v\in Z$, called \emph{twin vertices}, with $c(u)=c(v)\in \{i,j\}$, such that for every proper colouring $c'$ satisfying $c(y)\ne c'(y)\in \{i,j\}$ for some $y\in \{u,v\}$ and $c'(x)=c(x)$ for every $x\in Z\setminus\{y\}$, the sequence $Z$ (up to a permutation) is a zigzag in $(G,c')$.

The following two observations are used often (sometimes implicitly) throughout the proof.

\begin{observation}\label{obs:maxzig}
Let $c$ be a proper $t$-colouring of a graph $G$ with $\zig(G)=t$. Every maximum zigzag in $(G,c)$ is of size exactly $t$.
\end{observation}

\begin{proof}
    Since $\zig(G)=t$ we have $\zig(G,c)\ge t$, and since $c$ is a $t$-colouring we have $\zig(G,c)\le t$. 
\end{proof}

\begin{observation}\label{obs:verticesInPotentialZigzag}
Let $c$ be a proper $t$-colouring of a graph $G$ with $\zig(G)=t$ and $Z = z_1,z_2,z_3, \dots z_t$ be an $\{i,j\}$-potential zigzag of $(G,c)$. Then for all $k \in [t]\setminus\{i,j\}$, $c(z_k)=k$ and we have exactly one of the following:
$c(z_i) = c (z_j) = i$ or $c(z_i) = c(z_j) =j$. 
\end{observation}

\begin{lemma} \label{claim:potzig} 
Let $c$ be a proper $t$-colouring of a graph $G$ with $\zig(G)=t$ and $C$ be the union of some components of $G_{ij}$. Let $c'$  is the proper colouring of $G$ obtained from $c$ by switching colours $i, j$ in $C$. Let $Z$ be an $\{i,j\}$-potential zigzag in $(G,c)$. If $Z$ has exactly one of its vertices in $C$, then $Z$ (up to a permutation) is a maximum zigzag in $(G,c')$. Moreover, we have exactly one of the following, for any maximum zigzag $Z'$ of $(G,c')$.
\begin{itemize}
    \item $Z'$ has exactly one vertex in $C$ and $Z'$ is an $\{i,j\}$-potential zigzag $Z$ in $(G,c)$.
    \item $Z'$ (up to permutation) is a maximum zigzag in $(G,c)$.
\end{itemize}

\end{lemma}

\begin{proof}
Let $u,v$ be the pair of twin vertices in $Z$. Suppose $Z$ has exactly one vertex in $C$. The definition of $c'$ and $G_{ij}$ guarantees that either $u$ or $v$ is in $C$, say it is $u$. Thus we have $c'(u) \neq c(u)$ and $c'(v) = c(u)$. Moreover, since $u,v$ are the twin vertices in $Z$, we have $\{c'(u),c'(v)\} = \{i,j\}$. The definition of $\{i,j\}$-potential zigzag then gives us that $Z$ with some rearrangements of its elements is an maximum zigzag of $(G,c')$.
   
    Let $u'$ and $v'$ be the vertices of $Z'$ with colours $i$ and $j$, respectively, in the colouring $c'$. As $c'$ is obtained from $c$ by switching the colours of the vertices in $C$, each vertex in $Z'\setminus \{u',v'\}$ has the same colour in both the colouring $c$ and $c'$.  We have three cases depending on whether $u'$ or $v'$ is in $C$. If neither $u'$ nor $v'$ is a vertex of $C$. Then, 
    $Z'$ clearly is a maximum zigzag in $(G,c)$. If both $u'$ and $v'$ lie in $C$, then both $u'$ and $v'$ have switched their  colours. That is $c(u')=j$ and $c(v')=i$. Therefore, if $Z$ is the sequence obtained by interchanging the position of $u'$ and $v'$ in $Z'$, then $Z$ is a maximum zigzag in $(G,c)$.  
    
    Now we may assume that exactly one of the vertices among $u'$ and $v'$ is in $C$, without loss of generality, say $u'$. Then $u'$ has changed its colour but not $v'$. That is, $c(v')=c'(v')=j$ and $c'(u')=i$, $c(u')=j$, i.e., $Z'$ is a potential zigzag in $(G,c)$ with twin vertices $u',v'$.
\end{proof}

We can start now the proof of Theorem~\ref{theo:zigzag}. Let $G$ be a graph with $\chi(G)=\zig(G)=t$ and $c$ be a proper $t$-colouring of $G$ having a minimum number of maximum zigzags among $t$-colourings of $G$. Further let $Z = z_1, \dots z_t$ be a maximum zigzag in $(G,c)$. Note that by Observation~\ref{obs:maxzig}, we indeed know that $Z$ is of size $t$. 
Also note that since $c$ is a $t$-colouring we have 
$c(z_{i})=i$ for every $i \in [t]$. 
We choose as our set of terminals the vertices $\{z_i\in Z\mid i \mbox{ is odd}\}$. For each pair of terminals $z_i,z_j$, we will find  an odd length path $P_{ij}$ such that 
\begin{equation}\label{eq:colorescamino}
    \mbox{ for every $uv\in P_{ij}$ we have $\{c(u),c(v)\}\in \{\{i,j\},\{i,\ell_{ij}^i\},\{j,\ell_{ij}^j\}\}$,}
\end{equation}
where $\ell_{ij}^i, \ell_{ij}^j$ are even numbers chosen  in the way guaranteed by Lemma~\ref{lemma:existevennumber}. It is not hard to see that this choice of colours will guarantee that all paths are mutually edge-disjoint, so obtaining such paths will finish the proof.

Fix $i,j \in [t]$ different and odd. 
Suppose $z_i$ and $z_j$ are in the same component of $G_{ij}$, that is, $C_{ij}^{z_i}=C_{ij}^{z_j}$
. Then there is an odd $(z_i,z_j)$-path $P_{ij}$ in $G_{ij}$ such that $\{c(u),c(v)\} = \{i,j\}$, for every edge $uv$ in $P_{ij}$ (and thus $P_{ij}$ satisifes~\eqref{eq:colorescamino}). Therefore, we may assume that $C_{ij}^{z_i} \neq C_{ij}^{z_j}$.

\begin{claim}\label{claimseqofpaths}
For some integer $m\geq 2$, there exist pairwise disjoint subgraphs $\C_1, \ldots, \C_m$ of $G$ such that we have the following. Each of $\C_1, \ldots, \C_m$ is the union of some components of $G_{ij}$,  with $\C_1 = C_{ij}^{z_i}$ and $z_j$ is a vertex in $\C_m$. For all $1 \leq p \leq m-1$ we have: 

\qitem{I} The set $\mathcal{Z}_p := \{Z' \mid Z' \text{ is a zigzag in } (G,c_p) \text{ and } Z' \text{ is an } \{i,j\}\text{-potential zigzag in } (G,c)\}$ is  non-empty, where $c_p$ is the colouring of $G$ obtained from $c$ by switching colour $i$ and colour $j$ in $\C_1 \cup \C_2 \cup \ldots \cup \C_p$. Moreover we have:

\qitem{II} Let $a,b$ be the pair of twin vertices in an $\{i,j\}$-potential zigzag of (G,c) in $\mathcal{Z}_p$.
Then, exactly one of $a$ and $b$ is in $\C_p$ and the other is not in $\C_1 \cup \C_2 \cup \ldots \cup \C_p$.

\end{claim}

\begin{claimproof}
Set $\C_1 = C_{ij}^{z_i}$. We prove properties I and II by induction on $p$. For the base case $p=1$, we first show I, that is, $\mathcal{Z}_1$ is non-empty. As we have $C_{ij}^{z_i} \neq C_{ij}^{z_j}$, then 
in the colouring $c_1$ 
vertices $z_i$ and $z_j$ have the same colour, $j$. Hence $Z$, even with a permutation of its elements, is not a maximum zigzag in $(G,c_1)$. Now suppose for a contradiction that $\mathcal{Z}_1$ is an empty set, that is, each maximum zigzag in $(G,c_1)$ is also a maximum zigzag in $(G,c)$. Since $Z$ is not a maximum zigzag in $(G,c_1)$, we have that $(G,c_1)$ has fewer maximum zigzags than $(G,c)$. This contradicts the choice of $c$. Therefore, $\mathcal{Z}_1$ is a non-empty set. 

We now show property II for $p=1$. Let $Z'$ be any element of $\z_1$. Then $Z'$ is an $\{i,j\}$-potential in $(G,c)$. Let $a,b$ be the corresponding twin vertices. Then we claim that exactly one of the vertices among $a$ and $b$ is in $\C_1$. Suppose for a contradiction that either both $a$ and $b$ are in $\C_1$ or none of them. Given that $\{c_1(a),c_1(b)\} = \{i,j\}$ in both of these cases we obtain that $Z'$ is a maximum zigzag in $(G,c)$ if and only if $Z'$ (up to a permutation) is a maximum zigzag in $(G,c_1)$. This contradicts the assumption that $Z'$ is an element of $\z_1$. Thus, each element in $\z_1$ has exactly one twin vertex in $\C_1$. This verifies II for $p=1$, and thus the base case is ready. 
This base case is illustrated in Figure \ref{img::iteration2}.

\begin{figure}[H]
    \centering
\includegraphics[width=10cm]{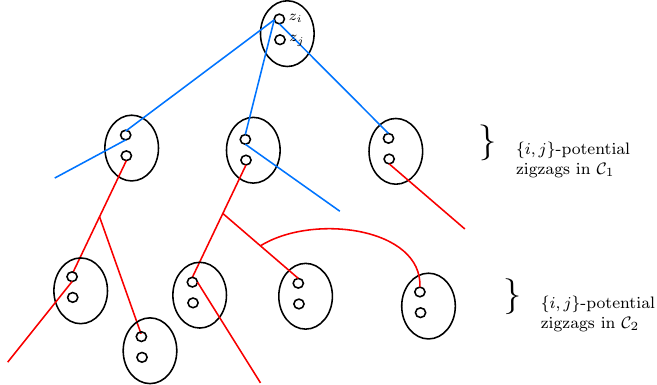}
    \caption{Base case of Claim \ref{claimseqofpaths}. }
    \label{img::iteration2}
\end{figure}

Suppose now that we have $m\ge p+1\ge 3$ and pairwise  disjoint subgraphs $\C_1, \C_2, \ldots, \C_{p}$ of $G_{ij}$ such that each of them is the union of components of $G_{ij}$ and each satisfies properties I and II. 
Set $A_q := \{a \mid a \in V(\C_q) \text{ is a twin of a potential zigzag in } \mathcal{Z}_q\}$ and $B_q := \{b \mid b \not\in V(\C_1 \cup \C_2 \cup \ldots \cup \C_q) \text{ is a twin of a potential zigzag in } \mathcal{Z}_q\}$, for all $1 \leq q \leq p$. Define $\C_{p+1}$ to be the union of the components in the set $\{C_{ij}^b \mid b \in B_p\}$, which is non-empty, since $\mathcal{Z}_p$ is not empty. 


Let us verify the claim for $p+1$. First, suppose for a contradiction that $\mathcal{Z}_{p+1}$ is an empty set, i.e., there is no maximum zigzag in $(G,c_{p+1})$ that is an $\{i,j\}$-potential zigzag in $(G,c)$.  Then, by Lemma \ref{claim:potzig}, each maximum zigzag in $(G,c_{p+1})$ is also a  maximum zigzag in $(G,c)$. This implies that $c_{p+1}$, which is also a $t$-colouring, has fewer maximum zigzags than $c$, contradicting the choice of $c$. Thus, $\z_{p+1}$ is an non-empty set.

Let $Z_{p+1}$ be any element in $\z_{p+1}$ and let $\{a,b\}$ be its twin vertices under $(G,c)$. We have $c(a) = c(b)$ and $c(a) \in \{i,j\}$. Since $Z_{p+1}$ is a maximum zigzag in $(G,c_{p+1})$ we also have $\{c_{p+1}(a),c_{p+1}(b)\} = \{i,j\}$. As a starting point, we claim that exactly one of the vertices among $a$ and $b$ is in $\C_{1} \cup \C_{2} \cup \ldots \cup \C_{p+1}$. Suppose for a contradiction that either both or none among $a$ and~$b$ lie in $\C_{1} \cup \C_{2} \cup \ldots \cup \C_{p+1}$. Then, $\{c_{p+1}(a),c_{p+1}(b)\} = \{c(a),c(b)\}$. This is a contradiction. Thus, exactly one of the vertices among $a$ and $b$ is in $\C_1 \cup \C_2 \ldots \cup \C_{p+1}$, say $a$.

We now arrive at showing that $a$ is in $\C_{p+1}$. Suppose for a contradiction that $a$ is a vertex in $\C_q$, for some $1 \leq q \leq p$. Then from the definition of $c_q$ and $b$ being a vertex not in $\C_1 \cup \C_2 \ldots \cup \C_{p+1}$, we know that $c_q(a) \neq c(a)$ and $c_q(b) = c(b)$. Thus, we have $\{c_q(a),c_q(b)\} = \{i,j\}$ and that $Z_{p+1}$ is also an element of $\z_q$. Hence, $b$ is a vertex in $\C_{q+1}$, a contradiction. Therefore, $a$ is a vertex in $\C_{p+1}$ and $b$ is not in $\C_1 \cup \C_2 \ldots \cup \C_{p+1}$. In other words, each $\{i,j\}$-potential zigzag in $\z_{p+1}$ has one twin vertex in $\C_{p+1}$ and one twin not  in $\C_1 \cup \C_2 \ldots \cup \C_{p+1}$. 

Since $G_{ij}$ has finite number of vertices, we will eventually have $z_j$ in $C_m$ for some value of $m\ge 2$. The result follows.
\end{claimproof}

\begin{figure}
    \centering
    \includegraphics[width= 9cm, height = 9cm]{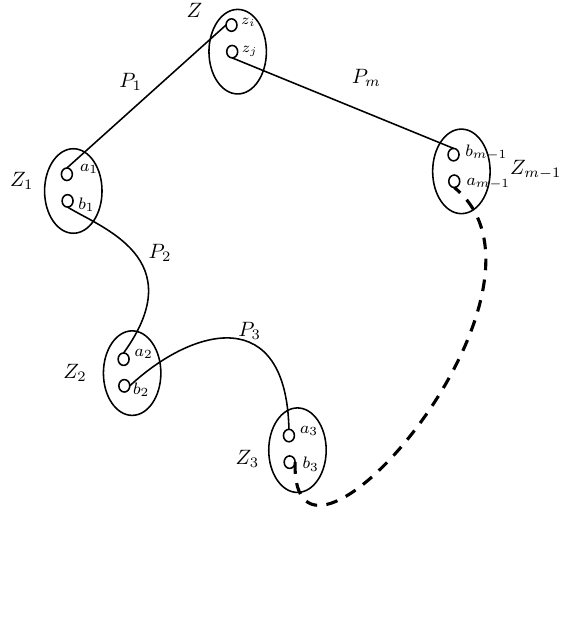}
    \caption{Paths $P_1,P_2,\dots, P_m$}
    \label{fig:Claim16}
\end{figure}

We now construct a sequence $P_1, P_2,\dots , P_{m}$ of pairwise (vertex-)disjoint paths in $G_{ij}$
and pick a sequence of $\{i,j\}$-potential zigzags $Z_1, Z_2, \dots , Z_{m-1}$ of $(G,c)$ with twin vertices $\{a_1, b_1\}, \,$ $ \{a_2, b_2\},\dots ,\{a_{m-1}, b_{m-1}\}$, respectively, 
such that
\begin{itemize}
    \item $P_1$  is a $(z_i,a_1)$-path,
    \item   $P_{m}$ is a $(b_{m-1},z_j)$-path, and  
    \item$P_q$ is a $(b_{q-1},a_q)$-path, for all $q\in \{2, \ldots, m-1\}$.
\end{itemize} 
See Figure~\ref{fig:Claim16}. By Claim~\ref{claimseqofpaths}, for every $1 \leq q \leq m-1$, each zigzag in $\z_q$ has one twin in $\C_q$ and the other twin  not in $\C_1 \cup \C_2 \cup \dots \cup \C_q$. As in the proof of that claim, let $A_q$ and $B_q$ denote $\{a \mid a \in V(\C_q) \text{ is a twin vertex of a zigzag in } \z_q\}$ and $\{b \mid b \not\in V(\C_1 \cup \C_2 \ldots \cup \C_q)$ is a twin vertex of a zigzag in $\z_q\}$, respectively, for all $1 \leq q \leq m-1$. Further let $b_{m-1}\in B_{m-1}$ be such that $z_j\in C_{ij}^{b_{m-1}}$. Choose $P_{m}$ as any $(b_{m-1},z_j)$-path in $C_{ij}^{b_{m-1}}$. Let $Z_{m-1}$ be a potential zigzag in $\mathcal{Z}_{m-1}$ with twin vertices $\{{a}_{m-1}, b_{m-1}\}$. Note that $a_{m-1}$ is a vertex of $A_{m-1}$, and that there exists a sequence of $\{i,j\}$-potential zigzags $Z_{m-2}, Z_{m-3}, \dots ,Z_{1}$ in $\z_{m-2}, \z_{m-3}, \ldots \z_1$, respectively, with twin vertices $\{a_{m-2},b_{m-2}\}, \{a_{m-3},b_{m-3}\}, \dots ,\{a_1,b_1\}$, respectively, where $a_{q} \in A_q$ and $b_q \in B_q$, for all $1 \leq q \leq m-2$. For $2\le q\le m-1 $ choose for $P_{q}$ some $(a_{q+1},b_{q})$-path in the subgraph $\C_{q}$
Finally, 
as $a_{1} \in \C_1 = C_{ij}^{z_i}$, we can choose a path $P_{1}$ in $\C_1$ joining $z_i$ and $a_{1}$. As $\C_1, \C_2, \ldots \C_{p+1}$ are pairwise disjoint subgraphs of $G_{ij}$, these paths are (vertex-)disjoint.

Now we will sew these paths together into a walk, out of which we will obtain our desired odd $(z_i,z_j)$-path. For every $1 \leq q \leq m-1$ let $v_q$ be the vertex in $Z_q$ of colour $\ell_{ij}^{c(a_q)}$, where $\ell_{ij}^i, \ell_{ij}^j$ are the even numbers given in Lemma~\ref{lemma:existevennumber}. Let $W$ be the walk obtained by connecting the following paths and vertices in the following order: $P_1-v_1-P_2-v_2-P_3-v_3- \ldots P_{m-1}-v_{m-1}-P_m$. 
By construction we have that for every edge $uv\in W$ we have $\{c(u),c(v)\}\in \{\{i,j\},\{i,\ell_{ij}^i\},\{j,\ell_{ij}^j\}\}$. Thus if we show that $W$ indeed contains an odd $(z_i,z_j)$-path, then this path will satisfy \eqref{eq:colorescamino}.

Define the colouring $\varphi\colon V(W)\rightarrow \{i,j\}$ as follows.

\begin{equation*}
\varphi(u) =
    \begin{cases}
        i, \text{ if } u \in \{v_1,v_2,\ldots v_{m-1}\} \text{ and } c(v_q)=\ell_{ij}^j\\
        j, \text{ if } u \in \{v_1,v_2,\ldots v_{m-1}\} \text{ and } c(v_q)=\ell_{ij}^i\\
        c(u), \text{ otherwise.} \\
    \end{cases}
\end{equation*}

Let us see that $\varphi$ is a proper colouring of $W$. Let $uv$ be an edge of $W$. If $uv$ is an edge in $P_q$, for some $1 \leq q \leq m$, then $\varphi(u) \neq \varphi(v)$, since $c$ is a proper colouring of $G$. So we may assume that $u=v_q$, for some $1 \leq q \leq m-1$ and, without loss of generality, $\varphi(v_q)=i$. Thus, $c(v_q) = \ell_{ij}^j$. The only neighbours of $v_q$ in $W$ are the pair of twin vertices $a_q$ and $b_q$. The definition of $W$ says that $c(a_q) = c(b_q) =j$, and thus $\varphi$ is a proper $2$-colouring of $H$.

Since $W$ is a walk from $z_i$ to $z_j$, $W$ contains a path $P_{ij}$ from $z_i$ to $z_j$. By definition, $\varphi$ assigns colour $i$ to $z_i$ and $j$ to $z_j$. Since $\varphi$ is a $2$-colouring, $P_{ij}$ is an odd path. The result follows.

\section{Odd immersions in Generalised Mycielski graphs} \label{sec:submicyelski} 

\noindent
In this section, we prove Theorem~\ref{theo:Mycielski}. Straightforward modification of this proof gives also Theorem~\ref{theo:Mycielski2}.

Let $G$ be a graph containing $K_t$ as a totally odd immersion with terminals $T = \{v_1,v_2,\dots, v_t\}$ and odd paths $P_{ij}$, for $1\le i<j\le t$. 
We will show that $\mu_m(G)$ has a totally odd immersion with set of terminals $\tau_1, \tau_2, \dots ,\tau_t,\tau_*$, a $(\tau_i,\tau_j)$-path $Q_{ij}$ for every  $1\le i<j \le t$, and a $(\tau_i,\tau_*)$-path $Q_{i}$ for every $1\le i\le t,$ all of which will be defined accordingly in each case. 

We first consider the easier case when there is no $P_{ij}$ of length 1. If $m$ is even we let $\tau_i=(v_i,1)$ for every $1\le i\le t$, and otherwise let $\tau_j=(v_i,0)$. Set $\tau_*=w$. For every pair  $1\le k<k'\le t$, let $a_1^{kk'},\dots , a_{2\ell}^{kk'}\in V(G)$ be such that $P_{kk'}=v_ka_1^{kk'},\dots , a_{2\ell}^{kk'}v_{k'}$, and take $Q_{kk'} := \tau_k,(a_1^{kk'},0), \dots ,(a_{2\ell}^{kk'},0),\tau_{k'}$. Notice that for any $k,k',s,s' \in \{1, \dots t\}$, $Q_{kk'}$ is edge disjoint with $Q_{ss'}$, whenever $\{k,k'\} \neq \{s,s'\}$. Moreover, since every path $P_{kk'}$ is of length at least three, we can choose, for every $1\le i\le t$, a vertex $x_i\in N_G(v_i)$ such that $x_i\ne x_j$ whenever $i\ne j$. For every $1 \leq i \leq t$ take $$Q_i = (v_i,1),(x_i,2),(v_i,3),\dots ,(x_i,m-2),(v_i,m-1),w, \mbox{ if }m\mbox{ is even, and }$$ $$Q_i = (v_i,0),(x_i,1),(v_i,2),\dots ,(v_i,m-2),(x_i,m-1),w, \mbox{ otherwise.}$$ These paths are mutually internally vertex-disjoint and also internally vertex-disjoint from all the  $Q_{ij}$ paths. Thus we have obtained the desired totally odd immersion.

Now we consider the case in which there exist a pair $i,j \in \{1,2, \ldots t\}$ such that $P_{ij}$ is an edge. Without loss of generality, we may assume that $j=1$. We set $\tau_i=(v_i,0)$ for every $1\le i\le t$ and $\tau_*=(v_1,1)$. For every pair of distinct $k,k'\in \{2,3, \ldots t\}$, let $a_1^{kk'},\dots , a_{2\ell}^{kk'}\in V(G)$ be such that $P_{kk'}=v_ka_1^{kk'},\dots , a_{2\ell}^{kk'}v_{k'}$, and take $Q_{kk'} := (v_k,0),(a_1^{kk'},0), \dots ,(a_{2\ell}^{kk'},0),(v_{k'},0)$. Notice that for any $k,k',s,s' \in \{2,3, \dots t\}$, $Q_{kk'}$ and $Q_{ss'}$ are mutually edge-disjoint whenever $\{k,k'\} \neq \{s,s'\}$.

For every $k\in \{2,\dots ,t\}$ we now join $(v_k,0)$ to both $(v_1,0)$ and $(v_1,1)$.  If $P_{1k}$ is an edge, then we take the path $Q_{1k}$ to be the edge, $(v_1,0)(v_k,0)$, and $Q_{k}$ to be the edge $(v_1,1)(v_k,0)$. Otherwise, we know $P_{1k}$ is of the form $v_1,b_1^{1k}, \dots, b_{2\ell'}^{1k},v_k$, for some integer $\ell'\ge 1$, and take $$Q_{1k} := (v_1,0),(b_1^{1k},1),(b_2^{1k},0),(b_3^{1k},1), \dots (b_{2\ell'-1}^{1k},1),(b_{2\ell'}^{1k},0),(v_k,0), \mbox{ and}$$ $$Q_{k} := \newline (v_1,1),(b_1^{1k},0),(b_2^{1k},1),(b_3^{1k},0), \dots ,(b_{2\ell'-1}^{1k},0),(b_{2\ell'}^{1k},1),(v_k,0).$$ Note that $(b_{x},y)$ is a vertex in $Q_{1k}$ if $x$ and $y$ are of same parity and $(b_{x},y)$ is a vertex in $Q_{k}$ if $x$ and $y$ are of different parity. Therefore, $Q_{1k}$ and $Q_{k}$ are internally vertex-disjoint. Also, if $(x,y)(x',y')$ is an edge of $Q_k$, then we have $y\ne y$, which guarantees that $Q_k$ is edge-disjoint from $Q_{ss'}$ for every pair of distinct $s,s'\in \{2,\dots , t\}$. To see that $Q_{1k}$ is edge-disjoint from every $Q_{ss'}$ we use a similar argument, plus the fact that if $(b,0)(v_k,0)\in E(Q_{1k})$ and $(b',0)(v_k,0)\in E(Q_{ss'})$, then we have $b\ne b'$ because $P_{1k}$ is edge-disjoint from $P_{ss'}$ and we have $bv_k\in P_{1k}, b'v_k\in P_{ss'}$.

Now we only need to find the path $Q_{1}$ which joins $(v_1,0)$ and $(v_1,1)$. Recall that $v_i$ is a neighbour of $v_1$. If $m$ is odd, then we take $Q_1 := (v_1,0),(v_i,1),(v_1,2), \dots ,(v_1,m-1),w,(v_i,m-1),(v_1,m-2),\dots   ,(v_1,1)$. Otherwise, we take  $Q_1 := (v_1,0),(v_i,1),(v_1,2), \ldots, (v_i,m-1),w,(v_1,m-1),(v_i,m-2),\dots   ,(v_1,1)$. Note that $Q_1$ is of length $2m-1$. Since the edge $v_1v_i$ does not belong to any path $P_{1s}$ with $s\ne i$, then $Q_1$ is internally vertex-disjoint from $Q_{1s}$. Moreover if $Q_{1i}$ has length 1, then we have $Q_{1i}=(v_1,0),(v_i,0)$, and otherwise its first edge is of the form $(v_1,0)(b_1,1)$ for $b_1\ne v_i$. It follows that $Q_1$ is internally vertex-disjoint from all other chosen paths. The result follows.

\section*{Acknowledgments}

\noindent
 H. Echeverría is supported by ANID BECAS/DOCTORADO NACIONAL 21231147.  
  A. Jim\'enez is  supported by  ANID/Fondecyt Regular 1220071 and ANID-MILENIO-NCN2024-103.
  A. Pastine is supported  by PIP CONICET 11220220100068CO, PICT 2020-00549, PICT 2020-04064,  PROICO 03-0723, and PROIPRO 03-2923.
S. Mishra is partially supported by ANID/Fondecyt Postdoctoral grant $3220618$ of Agencia Nacional de Investigati\'{o}n y Desarrollo (ANID), Chile.
  D.A. Quiroz is supported by ANID/Fondecyt Regular 1252197 and by MATH-AMSUD MATH230035.
  M. Yépez is supported by ANID BECAS/DOCTORADO NACIONAL 21231444.

\section*{Appendix}
\appendix
\section{$r\geq k=3$} 

\noindent
Recall that we want to find a collection of pairwise internally vertex-disjoint paths of the form $X_iB_{i,j}C_{i,j}X_j$ joining all pairs of vertices of the form
$X_i=\{k+1,\ldots,2k-1\}\cup\{2k+i\},$
with $0\leq i\leq r$.
For this case, we choose the paths similarly as in the case $r\ge k\ge 4$, but a little more carefully.

We first consider $0<i<j<r$. For $j-i=1$, we have four cases.

\begin{itemize}
\item If $i=1$ and $r=3$ 
take $B_{i,j}=B_{1,2}=\{1\}\cup\{6,8\}, C_{i,j}=C_{1,2}=\{2\}\cup \{7,9\}.$

\item If $i=1$ and $r=4$ set $B_{i,j}=B_{1,2}=\{1\}\cup\{6,8\},
    C_{i,j}=C_{1,2}=\{3\}\cup \{9,10\}.$

\item If $i=1$ and $r\ge 5$ take $B_{i,j}=B_{1,2}=\{1\}\cup\{6,8\}, C_{i,j}=C_{1,2}=\{3\}\cup \{7,10\}.$

\item otherwise, take $B_{i,j}=B_{i,i+1}=\{1\}\cup\{5+i,7+i\},
    C_{i,j}=C_{i,i+1}=\{3\}\cup \{4+i,8+i\}.$

\end{itemize}

For $j-i\geq 2$, we set $B_{i,j}=\{1\}\cup\{5+i,6+j\},
C_{i,j}=\{2\}\cup\{6+i,5+j\}.$

For the $(X_0, X_i)$-paths with $i\le r-1$, we have two cases, as follows.

 \begin{itemize}
     \item If $j\leq 2$, set $B_{0,j}=[1,2]\cup\{6+j\}, C_{0,j}=\{3\}\cup[7,5+j]\cup[7+j,9].$

\item If $3\leq j<r$, set
$B_{0,j}=[1,2]\cup\{6+j\},
C_{0,j}=\{3\}\cup\{5+j\}\cup\{7+j\}.$

 \end{itemize}
 
For the $(X_i, X_r)$-paths with $i\le r-1$, we also have two cases.
\begin{itemize}
    \item if $i= r-1$, set $B_{i,r}=[1]\cup\{4+r,6+r\}, 
C_{i,r}=[2,3]\cup\{5+r\}.$

\item if $i<r-1$, set
$B_{i,r}=[1]\cup\{5+i,6+r\},
C_{i,r}=[2,3]\cup\{6+i\}.$
\end{itemize}
Finally, set 
$B_{0,r}=[1,2]\cup\{6+r\}$ and
$C_{0,r}=\{3\}\cup[4+r,5+r].$

Sets $B_{ij}$ and $C_{ij}$, and the desired paths are uniquely determined as before.

\section{$k=2$}  The techniques used in the other cases where $r\ge k$ cannot be applied here since there are not enough vertices to make so many paths of length 3. Instead, we consider the following $r+2$  terminals:
$X_i=\{2i-1,2i\}$, for $i\in \{1,2,\ldots \lfloor \frac r2 \rfloor +2\},$ and 
$Y_i=\{2i,2i+1\}$, for $i\in \{2,\ldots,\lceil \frac r2 \rceil+1\}$. The set of all $X_i$ forms a clique; so does the set of all $Y_i$. Further, each $Y_i$ is adjacent to every $X_j$ except for $X_i$ and $X_{i+1}$. To join these non-adjacent terminals,
when $r$ is even, we consider the paths $$P_i=Y_i,\{1,2i-1\},\{2,2i+1\},X_i$$ $$Q_i=Y_i,\{1,2i+2\},\{2,2i\},X_{i+1}$$
for all $i\in \{2,\ldots,\lceil \frac r2 \rceil+1\}$.
Then, every path $P_i$ is vertex disjoint from all the paths $Q_j$. Finally, for $i\neq j$, the path $P_i$ is disjoint from  $P_j$, and the same holds for $Q_i$ and $Q_j $.

The case when $r$ is odd works in  the same way, but we note that when $i=\lceil \frac r2 \rceil+1$ there is no $X_{i+1}$, and thus the path $Q_i$ is not needed.

\end{document}